\documentstyle[12pt]{article}
\newcommand{\sect}[1]{\section{#1}\setcounter{equation}{0}}

\font\mbn=msbm10 scaled \magstep1
\font\mbs=msbm7 scaled \magstep1
\font\mbss=msbm5 scaled \magstep1
\newfam\mbff
\textfont\mbff=\mbn
\scriptfont\mbff=\mbs
\scriptscriptfont\mbff=\mbss
\def\mbf{\fam\mbff}
\def\Re{{\mbf R}}

\def\Co{{\mbf C}}

\def\Di{{\mbf D}}
\newtheorem{Th}{Theorem}[section]
\newtheorem{Lm}[Th]{Lemma}
\newtheorem{C}[Th]{Corollary}
\newtheorem{D}[Th]{Definition}
\newtheorem{Proposition}[Th]{Proposition}
\newtheorem{R}[Th]{Remark}

\author{Alexander Brudnyi\\
Department of Mathematics\\
Ben Gurion University of the Negev, Beer-Sheva\\
Israel}
\title{Traces of Subharmonic Functions to Fractal Sets}
\date{June 23, 1999}
\begin{document}
\maketitle
\begin{abstract}
We study traces of a class of subharmonic functions to Ahlfors regular
subsets of $\Co^{n}$. In particular, we establish for the traces a generalized
BMO-property and the reverse H\"{o}lder inequality.
\end{abstract}
\sect{\hspace*{-1em}. Introduction.}
{\bf 1.1.} 
A compact subset $K\subset\Re^{n}$ is said to be (Ahlfors) {\em $d$-regular}
if there is a positive number $a$ such that for any $x\in K,\ \ 0<t\leq
{\rm diam}(K)$
\begin{equation}\label{alf1}
{\cal H}^{d}(K\cap\Di(x,t))\leq at^{d}.
\end{equation}
Here ${\cal H}^{d}(\omega)$ denotes the $d$-Hausdorff measure of $\omega$.
This class will be denoted by ${\cal A}(d,a)$.

A compact subset $K\in {\cal A}(d,a)$ is said to be a {\em $d$-set}
if there is a positive number $b$ such that for any $x\in K,\ \ 0<t\leq
{\rm diam}(K)$
\begin{equation}\label{alf2}
bt^{d}\leq {\cal H}^{d}(K\cap\Di(x,t)).
\end{equation}
Denote this class by ${\cal A}(d,a,b)$.

The purpose of this paper is to study traces
of subharmonic functions to $d$-sets. Let us recall that the class of $d$-sets,
in particular, contains Lipschitz $d$-manifolds (with $d$ integer), Cantor
type sets and self-similar sets (with arbitrary $d$), see, e.g., [JW], p. 29
and [M], sect. 4.13.

Let us, first, formulate our results in $\Co$. In the sequel we
denote $\Di_{s}:=\{z\in\Co:\ |z|<s\}$ and $\Di(x,t):=\{z\in\Co\ :\ |z-x|<t\}$.

Assume that $f$ is a subharmonic in $\Di_{1}$ function satisfying
\begin{equation}\label{m1m2}
\sup_{\Di_{1}}f\leq M_{1}\ \ {\rm and}
\ \ \ \sup_{\Di_{r}}f\geq M_{2}\ \ \ \ \ \ (r<1).
\end{equation}
\begin{Th}\label{rem}
Let $\omega\subset\Di(x,t)$ be a compact
set of ${\cal A}(d,a)$ satisfying ${\cal H}^{d}(\omega)\geq\epsilon>0$.
Assume that $\Di(x,t/r)\subset\Di_{r}$.
Then there is a constant $c=c(r)>0$ such that inequality
$$
\sup_{\Di(x,t)}f\leq\sup_{\omega}f+(M_{1}-M_{2})c\log
\frac{4eta^{1/d}}{r(d\epsilon)^{1/d}}
$$
holds for any subharmonic $f$ satisfying (\ref{m1m2}).
\end{Th}
We use Theorem \ref{rem} to establish a generalized BMO-property and the
reverse H\"{o}lder inequality for traces of
subharmonic functions to $d$-sets. Let us recall
\begin{D}\label{bmo}
Let $X$ be a complete metric space equipped with a regular Borel
measure $\mu$. A locally integrable on $X$ function $f$ belongs
to $BMO(X,\mu)$ if
$$
|f|_{*}:=\sup\left\{\frac{1}{\mu(B)}\int_{B}|f-f_{B}|d\mu\right\}<\infty ;
$$
here supremum is taken over all metric balls $B\subset X$ and
$f_{B}=\frac{1}{\mu(B)}\int_{B}fd\mu$.
\end{D}
In order to formulate the next two results consider a compact $d$-set
$K\subset\Co$. Assume that $f$ is a subharmonic function defined in an open
neighbourhood of $K$.

\begin{Th}\label{tbmo}
Restriction $f|_{K}$ belongs to $BMO(K,{\cal H}^{d})$.
\end{Th}
\begin{Th}\label{reve}
For any $K_{x,t}:=\Di(x,t)\cap K$,
$x\in K$, $t>0$, and any $1\leq p\leq\infty$ the inequality
\begin{equation}\label{rev}
\left(\frac{1}{{\cal H}^{d}(K_{x,t})}\int_{K_{x,t}}e^{pf}d{\cal H}^{d}
\right)^{1/p}\leq C(K,f,d)\frac{1}{{\cal H}^{d}(K_{x,t})}\int_{K_{x,t}}e^{f}
d{\cal H}^{d}
\end{equation}
holds.
\end{Th}

{\bf 1.2.} In this section we consider multidimensional generalizations
of the previous results for plurisubharmonic ({\em psh}) functions in
$\Co^{n}$, $n\geq 2$. Note that it is impossible to obtain in this situation
the results similar to those in $\Co$. Indeed, let $N:=\{z\in\Co^{n}\ :\
f(z)=0\}\neq\emptyset$ for a nontrivial holomorphic function $f$. Let
$K=K_{1}\cup K_{2}$ be a connected set where $K_{1},\ K_{2}$ are compact
$(2n-2)$-Lipschitz manifolds such that
$K_{1}\subset N$ and $K_{2}\not\subset N$. Then the
psh function $\log |f|$ equals $-\infty$ on $K_{1}$ and $\sup_{K}\log |f|>
-\infty$. Hence the results similar to Theorems \ref{rem}, \ref{tbmo} and
\ref{reve} are not true in this case.

So we restrict ourselves to a special class of psh
functions whose local behavior is similar to that of subharmonic functions in
$\Co$.

Consider a family $F=\{f_{1},...,f_{n}\}$ of holomorphic functions
defined in Euclidean ball $B_{r}\subset\Co^{n}$ of radius $r>1/2$ 
centered at $0$. We will prove the results for
psh function $u=\log |F|$ where
$|F|^{2}:=|f_{1}|^{2}+...+|f_{n}|^{2}$. Assume that 
\begin{equation}\label{plur}
\sup_{B_{r}}u\leq 0 \ \ \ {\rm and}\ \ \ \inf_{\partial B_{1/2}}u\geq -M
\end{equation}
for some $M>0$. The latter condition means that the system $F=0$ has only
discrete zeros in a neigbourhood of the closed ball $\overline{B_{1/2}}$.
\begin{D}\label{ell}
Let $h$ be a nonnegative analytic function with a finite number of zeros
defined on an open set $U\subset\Re^{N}$.
A zero $x$ of
$h$ is said to be elliptic if the Taylor expansion of $f$ at $x$ has the
following form
$$
h(x+t\omega)=t^{d}f(\omega)+o(t^{d})
$$
with
$$
\inf_{S^{N-1}}f>0\ .
$$
Here $t>0$ and $\omega$ belongs to the unit sphere $S^{N-1}$ of $\Re^{N}$.
\end{D}
Assume that $F$ satisfies condition (\ref{plur})  and all zeros of
$|F|^{2}$ are elliptic. Let $k$ be the number of zeros of $|F|^{2}$ in
$B_{1/2}$ counting with their multiplicities.
(By multiplicity we mean the degree of map $F$ at zero.) 
\begin{Th}\label{mcol1}
Let $\omega\subset B(x,t)$  be a compact set of 
${\cal A}(d,a)$ and ${\cal H}^{d}(\omega)\geq\epsilon>0$.
Assume that $B(x,4r^{2}t)\subset B_{1/2}$. Then
there is a constant $c=c(r,F)>0$ such that
$$
\sup_{B(x,t)}\log |F|\leq\sup_{\omega}\log |F|+
k\log\frac{16erta^{1/d}}{c(d\epsilon )^{1/d}}\ .
$$
\end{Th}
The following theorem gives the results similar to those of
Theorems \ref{tbmo} and \ref{reve}.
\begin{Th}\label{mcol2}
Assume that a compact set $K\subset\Co^{n}$ belongs to ${\cal A}(d,a,b)$.
Then under assumptions of Theorem \ref{mcol1}, $|F|$ satisfies the reverse
H\"{o}lder inequality in each ball $B(x,t)\cap K$ and
$\log |F|\in BMO(K,{\cal H}^{s})$.
\end{Th}
\sect{\hspace*{-1em}. Abstract Version of Cartan's Lemma.}
Our proofs are based on estimates for psh
functions which generalize well-known Cartan's Lemma for polynomials
(see [Ca]).
We need a version of the generalized Cartan's Lemma proved by
Gorin (see [GK]).

Let $X$ be a complete metric space and let $\mu$ be
a finite Borel measure on $X$. We consider a continuous,
strictly increasing,
nonnegative function $\phi$ on $[0,+\infty [,\ \phi(0)=0,\ 
\lim_{x\to\infty}\phi(x)>\mu(X)$.  The function $\phi$ will be called a
${\em majorant}$.

For each point $x\in X$ we set $\tau(x)=\sup\{t:\ \mu(B(x,t))\geq\phi(t)\}$,
where $B(x,t)$ is the closed ball in $X$ with center $x$ and radius $t$.
It is easy to see that $\mu(B(x,\tau(x))=\phi(\tau(x))$  and
$\sup_{x}\tau(x)\leq\phi^{-1}(\mu(X))<\infty$.

A point $x\in X$ is said to be {\em regular} 
(with respect to $\mu$ and $\phi$) if $\tau(x)=0$, i.e., 
$\mu(B(x,t))<\phi(t)$ for all $t>0$. The next result shows that the set of 
regular points is sufficiently large for an arbitrary majorant $\phi$.
\begin{Lm}\label{Go}(Gorin)
Let $0<\gamma<1/2$. There exists a sequence of balls $B_{k}=B(x_{k},t_{k})$,
$k=1,2,...,$ which collectively cover all the irregular points and which are
such that
$\sum_{k\geq 1}\phi(\gamma t_{k})\leq\mu(X)$ (i.e., $t_{k}\to 0$).
\end{Lm}
For the sake of completeness we present Gorin's proof of the lemma.\\
{\bf Proof.}
Let $0<\alpha<1$, $\beta>2$ but $\gamma<\alpha/\beta$. We set 
$B_{0}=\emptyset$ and assume that the balls $B_{0},...,B_{k-1}$ have been
constructed. If $\tau_{k}=\sup\{\tau(x):\ x\not\in B_{0}\cup...\cup B_{k-1}\}$,
then there exists a point $x_{k}\not\in B_{0}\cup...\cup B_{k-1}$, such that
$\tau(x_{k})\geq\alpha\tau_{k}$. We set $t_{k}=\beta\tau_{k}$ and
$B_{k}=B(x_{k},t_{k})$. Clearly, the sequence $\tau_{k}$ (and thus also
$t_{k}$) does not increase. The balls $B(x_{k},\tau_{k})$ are pairwise 
disjoint. Indeed, if $l>k$, then $x_{l}\not\in B_{k}$, i.e., the distance
between $x_{l}$ and $x_{k}$ is greater than $\beta\tau_{k}>2\tau_{k}\geq
\tau_{k}+\tau_{l}$. Then,
$$
\sum_{k=1}^{\infty}\phi(\gamma t_{k})\leq\sum_{k=1}^{\infty}
\phi(\alpha\tau_{k})\leq\sum_{k=1}^{\infty}\phi(\tau(x_{k}))=
\sum_{k=1}^{\infty}\mu(B(x_{k},\tau_{k}))\leq\mu(X)\ ;
$$
consequently, $\tau_{k}\to 0$, i.e., for each point $x$, not belonging to the
union of the balls $B_{k}$, $\tau(x)=0$, $x$ is a regular point. In addition,
$t_{k}=\beta\tau_{k}\to 0$.\ \ \ \ \ $\Box$
\begin{R}\label{comp}
{\rm If $X$ is a locally compact metric space then one can take $\gamma=1/2$ 
(for similar arguments see, e.g., [L], Th. 11.2.3).}
\end{R}
We now apply Lemma \ref{Go} to obtain estimates for logarithmic potentials
of measures.

Assume that $X$
is a locally compact metric space with metric $d(.,.)$.
\begin{Th}\label{Cart}
Let
$$
u(z)=\int_{X}\log d(x,\xi)d\mu(\xi)
$$
where $\mu$ is a Borel measure, $\mu(X)=k<\infty$. 

Given $H>0, d>0$ there exists a system of metric balls such that
\begin{equation}\label{condi}
\sum r_{j}^{d}\leq\frac{(2H)^{d}}{d}
\end{equation}
where $r_{j}$ are radii of these balls, and
$$
u(z)\geq k\log\frac{H}{e}
$$
everywhere outside these balls.
\end{Th}
{\bf Proof.}
Let $\phi(t)=(pt)^{d}$ be a majorant with $p=\frac{(kd)^{1/d}}{H}$.
We cover all irregular points of
$X$ by balls according to Gorin's Lemma \ref{Go} and Remark \ref{comp}. 
It remains to estimate the potential $u$ outside of these balls, i.e., at
any regular point $z$. Let $n(t;z)=\mu(\{\xi:\ d(z,\xi)\leq t\})$.
Clearly, for any $N\geq\max\{1, H\}$
$$
u(z)\geq\int_{d(z,\xi)\leq N}\log d(z,\xi)d\mu(\xi)=\int_{0}^{N}\log t\
dn(t;z)=
n(t;z)\log t|_{0}^{N}-\int_{0}^{N}\frac{n(t;z)}{t}dt .
$$
Since $n(t;z)<(pt)^{d}$, we then have
$$
u(z)\geq n(N;z)\log N-\int_{0}^{N}\frac{n(t;z)}{t}dt.
$$
In addition, $n(t;z)\leq n(N;z)$ for $t\leq N$. Therefore,
$$
\begin{array}{c}
\displaystyle
u(z)\geq n(N;z)\log N-\int_{0}^{H}\frac{(pt)^{d}}{t}dt-
\int_{H}^{N}\frac{n(N;z)}{t}dt=\\
\displaystyle
n(N;z)\log N-\frac{(pH)^{d}}{d}-n(N;z)\log N+n(N;z)\log H=
-k+n(N;z)\log H
\end{array}
$$
Letting here $N\to\infty$ and taking into account
that $\lim_{N\to\infty}n(N;z)=k$ we obtain the required result.\ \ \ \ \ 
$\Box$
\sect{\hspace*{-1em}. Proofs of Theorems \ref{rem}, \ref{tbmo} and
\ref{reve}.}
{\bf Proof of Theorem \ref{rem}.}
We begin with
\begin{Proposition}\label{C1}
Let $u$ be a nonpositive subharmonic function on $\Di_{1}$ satisfying
$$
\sup_{\Di_{r}}u\geq -1\ \ \ \ \ {\rm for\ some}\ r<1.
$$
Then for any $H>0, d>0$ there is a set of disks such that 
\begin{equation}\label{e1}
\sum r_{j}^{d}\leq\frac{(2H)^{d}}{d},
\end{equation}
where $r_{j}$ are radii of these disks, and
$$
u(z)\geq c\log\frac{H}{e}
$$
outside these disks in $\Di_{r}$. Here $c=c(r)>0$ depends on $r$ only.
\end{Proposition}
{\bf Proof.}
Let $\kappa$ be a nonnegative radial $C^{\infty}$-function on $\Co$ 
satisfying
\begin{equation}\label{kern}
\int\int_{\Co}\kappa (x,y)dxdy=1\ \ \ and\ \ \
supp(\kappa)\subset \Di_{1}\ \ \ \ \ (z=x+iy).
\end{equation}
Let $u_{k}$ denote the function defined on $\Di_{1-1/k}$ by
\begin{equation}\label{approxim}
u_{k}(w):=\int\int_{\Co}\kappa (x,y)u(w-z/k)dxdy.\ \
\end{equation}
It is well known, see, e.g., [K], Theorem 2.9.2, that
$u_{k}$ is subharmonic on $\Di_{1-1/k}$ of the class $C^{\infty}$
and that
$u_{k}(w)$ monotonically decreases and tends to $u(w)$ for each
$w\in\Di_{1}$ as $k\to\infty$. Let $K:=\{z\in\Di_{1}\ :\
\frac{1+r}{2}\leq |z|\leq\frac{3+r}{4}\}$ be an annulus in $\Di_{1}$ and
$k\geq k_{0}=[\frac{8}{1-r}]+1$.
We are based on the following result (see, e.g., [Br], Lemma 2.3).

There are a constant $A=A(r)>0$ and numbers $t_{k}$, $k\geq k_{0}$, satisfying
$\frac{1+r}{2}\leq t_{k}\leq\frac{3+r}{4}$ such that
$u_{k}(z)\geq -A$ for any $z\in\Co,\ |z|=t_{k}$.

Then we can construct functions $f_{k}$ subharmonic on $\Co$ by
\[
f_{k}(z):=
\left\{
\begin{array}{cc}
\displaystyle
u_{k}(z)&(z\in\Di_{t_{k}});\\
\displaystyle
\max\left\{u_{k}(z),\frac{-2A\log|z|}{\log t_{k}}\right\}&
(z\in\Di_{1}\setminus\Di_{t_{k}});\\
\displaystyle
\frac{-2A\log|z|}{\log t_{k}}&(z\in\Co\setminus\Di_{1}).
\end{array}
\right.
\]
Without loss of generality we may assume that
$t_{k}\to t\in [\frac{1+r}{2},\frac{3+r}{4}]$ as $k\to\infty$.
Finally, define
$$
f(z)=(\overline{\lim_{k\to\infty}}f_{k}(z))^{*},
$$
where $g^{*}$ denotes upper semicontinuous regularization of $g$. Then
$f$ is subharmonic in $\Co$ satisfying
$$
f(z)=u(z)\ \ \ \ (z\in\Di_{1})\ \ \ \ {\rm and}\ \ \ \
f(z)=\frac{-2A\log|z|}{\log t}\ \ \ \ \ (z\in\Co\setminus\overline{\Di_{1}}).
$$
Consider now $\mu=\Delta f$. Then $\mu$ is a finite
Borel measure on $\Co$ supported in $\overline{\Di_{1}}$.
According to F. Riesz's theorem (see, e.g., [HK], Th. 3.9)
$$
\tilde f(z):=\frac{1}{2\pi}\int\int_{\Co}\log|z-\xi|d\mu(\xi)
$$
is subharmonic in $\Co$ and satisfies
$\Delta\tilde f=\Delta f=\mu$.
Thus $h=\tilde f-f$ is a real-valued harmonic in $\Co$ function.
Moreover, $h$ goes to infinity as
$\left(\frac{\mu(\Co)}{2\pi}-\frac{-2A}{\log t}\right)\log |z|$.
This immediately implies
$h=0$ and $\frac{\mu(\Co)}{2\pi}=\frac{-2A}{\log t}$.
Now according to Theorem \ref{Cart} applied to $f(=\tilde f)$,
for any $H>0,d>0$ there is a system of disks with radii $r_{j}$ satisfying
$\sum r_{j}^{d}\leq\frac{(2H)^{d}}{d}$
such that
$$
f\geq\frac{-2A}{\log t}\log\frac{H}{e}\geq\frac{-2A}{\log r}\log\frac{H}{e}
$$
outside these disks. It remains to set $c=\frac{-2A}{\log r}$.

The proof of the proposition is complete.\ \ \ \ \ $\Box$

Assume now that $f$ is subharmonic and satisfies (\ref{m1m2}). Then
by the main theorem in [Br] there is a constant $C=C(r)>0$ such
that the inequality
$$
\sup_{\Di(x,t/r)}f\leq C(M_{1}-M_{2})+\sup_{\Di(x,t)}f
$$
holds for any pair of disks $\Di(x,t)\subset\Di(x,t/r)\subset\Di_{r}$.	
Applying inequality of Proposition \ref{C1} to the function 
$$
u(z)=\frac{f(tz/r)-\sup_{\Di(x,t/r)}f}{C(M_{1}-M_{2})}\ \ \ \ \ 
(z\in\Di_{1})
$$ 
and then going back to $f$ we obtain the following
\begin{Proposition}\label{local}
There is a constant $c=c(r)>0$ such that for any disk $\Di(x,t)$ satisfying
$\Di(x,t)\subset\Di(x,t/r)\subset\Di_{r}$ and any $H>0, d>0$ there is a system
of disks such that
$$
\sum r_{j}^{d}\leq\frac{(2tH/r)^{d}}{d},
$$
where $r_{j}$ are radii of these disks, and
$$
f(z)\geq\sup_{\Di(x,t)}f+c(M_{1}-M_{2})\log\frac{H}{e}
$$
outside these disks in $\Di(x,t)$.
\end{Proposition}
\begin{R}\label{lev}
{\rm A particular case of Proposition \ref{local} for functions $u=\log|f|$
with holomorphic $f$ and for $d=1$ was proved in [L].}
\end{R}
We proceed to the proof of Theorem \ref{rem}.
First we show that $\omega$ can not be
covered by a system of disks such that
\begin{equation}\label{contra}
\sum r_{j}^{d}\leq\frac{(1-1/n)\epsilon}{2^{d}a}\ \ \ \ (n\geq 1)
\end{equation}
where $r_{j}$ are radii of these disks.
Assume to the contrary that there exists a system of disks
$\{\Di(x_{j},r_{j})\}$ whose radii satisfy (\ref{contra}) which covers
$\omega$. For any $x_{j}$ choose $y_{j}\in\omega$ so that
$|x_{j}-y_{j}|\leq r_{j}$. Then the system of disks $\{\Di(y_{j},2r_{j})\}$
also covers $\omega$. Since $\omega\in {\cal A}(d,a)$, we obtain
inequality
$$
{\cal H}^{d}(\omega)\leq\sum {\cal H}^{d}(\omega\cap\Di(y_{j},2r_{j}))\leq
2^{d}a\sum r_{j}^{d}<\epsilon
$$
which contradicts to ${\cal H}^{d}(\omega)\geq\epsilon $.

We now apply Proposition \ref{local} with 
$H_{n}=\frac{(d(1-1/n)\epsilon)^{1/d}r}{4ta^{1/d}}$. Since any system
of disks with $\sum r_{j}^{d}\leq\frac{(2tH_{n}/r)^{d}}{d}$ can not
cover $\omega$,  Proposition \ref{local} implies that there is a point
$x_{n}\in\omega$ such that
$$
\sup_{\omega}f\geq f(x_{n})\geq
\sup_{\Di(x,t)}f+c(M_{1}-M_{2})\log\frac{H_{n}}{e}
$$  
Letting $n\to\infty$ we get the required inequality.

Theorem \ref{rem} is proved. \ \ \ \ \ $\Box$

Our next result shows that $d$-regularity is a necessary condition
for the set to satisfy the inequality of Theorem \ref{rem}.
\begin{Proposition}\label{sharp}
Let $K\subset\Di_{1/2}$ be a compact set with ${\cal H}^{d}(K)<\infty$. Assume 
that the inequality
$$
\sup_{\Di(x,t)}f\leq\sup_{\omega}f+L+C\log
\frac{t}{\epsilon^{1/d}}
$$
holds for any $\omega\subset K\cap\Di(x,t)\subset\Di(x,3t/2)\subset\Di_{2/3}$,
$x\in K$, with ${\cal H}^{d}(\omega)=\epsilon$ and any $f$ subharmonic in
$\Di_{1}$ satisfying (\ref{m1m2}) with $r=2/3$ and some $M_{1},M_{2}$.
Here $L$ and $C>0$ depend on $K,d,M_{1},M_{2}$.
Then $K\in {\cal A}(d,c)$ for some $c>0$.
\end{Proposition}
{\bf Proof.}
For any $f$, $\omega$, $t\leq 1/9$ satisfying assumptions of the proposition
the inequality
$$
-C\log\frac{t}{\epsilon^{1/d}}\leq
\sup_{\Di(x,t)}f-\sup_{\omega}f-C\log
\frac{t}{\epsilon^{1/d}}\leq L<\infty
$$
holds. For a point $x\in K$ we set
$f_{x}(z)=\log|z-x|$ and $\epsilon_{t}:=
{\cal H}^{d}(\Di(x,t)\cap K)$. Clearly the family $\{f_{x}\}$ satisfies
inequality (\ref{m1m2}) with $r=2/3$, $M_{1}=3/2$ and $M_{2}=1/6$.
Then the above inequality applied to $f_{x}$ gets
$$
L\geq -C\log\frac{t}{\epsilon_{t}^{1/d}},
$$
that is equivalent to $\epsilon_{t}\leq \widetilde Lt^{d}$ for
$\widetilde L=e^{\frac{dL}{C}}$. So we checked the definition of $d$-regularity
for $t\leq 1/9$. For $t>1/9$ the inequality is obvious.\ \ \ \ \ $\Box$

Assume that $f$ satisfies (\ref{m1m2}) and $K\subset\Di_{r}$ is a compact
from ${\cal A}(d,a)$.
For a pair $\Di(x,t)\subset\Di(x,t/r)(\subset\Di_{r})$ we set
$K_{x,t}:=\Di(x,t)\cap K$ and $f_{x,t}=\sup_{\Di(x,t)}f$. Further, set 
$f'=f_{x,t}-f$.
In the proofs of Theorem \ref{tbmo} and  \ref{reve} we use 
\begin{Lm}\label{dis}
Let $D_{f'}(\lambda)={\cal H}^{d}\{y\in K_{x,t}\ :\ f'(y)\geq\lambda\}$ be
the distribution function of $f'$. Then 
\begin{equation}\label{new}
D_{f'}(\lambda)\leq\frac{(4et)^{d}a}{r^{d}d}e^{-\lambda d/(c(M_{1}-M_{2}))}.
\end{equation}
\end{Lm}
{\bf Proof.} The proof follows straightforwardly from the inequality 
of Theorem \ref{rem} where we choose $\omega:=D_{f'}(\lambda)$. We leave the
details to the reader.\ \ \ \ \ $\Box$

{\bf Proof of Theorem \ref{tbmo}.}
First, we prove a local version of the theorem. 
Assume that $K\subset\Di_{r}$ is a compact
from ${\cal A}(d,a,b)$ and $f$, $\Di(x,t)$ satisfy conditions of Lemma \ref{dis}.
From inequality \ref{new} it follows
\begin{equation}\label{estbmo}
\begin{array}{c}
\displaystyle
\frac{1}{{\cal H}^{d}(K_{x,t})}\int_{K_{x,t}}f'd{\cal H}^{d}\leq
\frac{1}{{\cal H}^{d}(K_{x,t})}\int_{0}^{\infty}D_{f'}(x)dx\leq
\frac{1}{bt^{d}}\frac{c(M_{1}-M_{2})}{d}\frac{(4et)^{d}a}{r^{d}d}=\\
\\
\displaystyle
\frac{ca(4e)^{d}(M_{1}-M_{2})}{br^{d}d^{2}}\ .
\end{array}
\end{equation}
Now we have
$$
\begin{array}{c}
\displaystyle
\frac{1}{{\cal H}^{d}(K_{x,t})}\int_{K_{x,t}}|f-f_{K_{x,t}}|d{\cal H}^{d}\leq
\frac{1}{{\cal H}^{d}(K_{x,t})}\int_{K_{x,t}}|(f-f_{x,t})-
(f-f_{x,t})_{K_{x,t}}|d{\cal H}^{d}\leq \\
\\
\displaystyle
\frac{2}{{\cal H}^{d}(K_{x,t})}\int_{K_{x,t}}f'd{\cal H}^{d}
\leq\frac{2ca(4e)^{d}(M_{1}-M_{2})}{br^{d}d^{2}}\ .
\end{array}
$$
This gives the estimate of the BMO-norm in each ball $K(x,t)=\Di(x,t)\cap K$
with $\Di(x,t)\subset\Di(x,t/r)(\subset\Di_{r})$.
In the general case,
we cover $K$ by a finite number of open disks $\Di(x_{i},R)$,
$i=1,...,N$ such that $f$ is defined in the union of these disks, the set
$\cup_{i=1}^{N}\Di(x_{i},R/2)$ also covers $K$  and any disk of radius 
$\leq R/4$ centered at a point of $K$ belongs to one of $\Di(x_{i},R/2)$.
Then the estimate of the 
BMO-norm in any $\Di(x,t)\cap K$, $x\in K$, $t\leq R/4$, follows
from Theorem \ref{rem} and inequality (\ref{estbmo}). To estimate
BMO-norms for $\Di(x,t)\cap K$ with $t\geq R/4$ we write   
$$
\frac{1}{{\cal H}^{d}(K_{x,t})}\int_{K_{x,t}}|f-f_{K_{x,t}}|d{\cal H}^{d}\leq
\frac{4^{d}}{bR^{d}}\int_{K_{x,t}}2|f|d{\cal H}^{d}<C\int_{K}|f|d{\cal H}^{d}.
$$
To complete the proof note that (\ref{new}) implies 
$\int_{K}|f|d{\cal H}^{d}<\infty$.\ \ \ \ \ $\Box$

We now formulate another corollary of Theorem \ref{rem}.
\begin{C}\label{ancol}
Assume that a subharmonic function $f$ defined on $\Co$ satisfies
$$
f(z)\leq c'+\log(1+|z|)\ \ \ \ \ (z\in\Co)\ 
$$
for some $c'\in\Re$.
Assume also that $S\in {\cal A}(d,a,b)$. Then $f|_{S}\in BMO(S,{\cal H}^{d})$
and the BMO norm $|f|_{S}|_{*}\leq\frac{\tilde ca}{bd^{2}}$ with an absolute
constant $\tilde c$.
\end{C}
{\bf Proof.}
For functions $f$ satisfying conditions of the corollary the 
Bernstein-Walsh inequality 
\begin{equation}\label{bernwal}
\sup_{\Di(x,qt)}f\leq\log q+\sup_{\Di(x,t)}f
\end{equation}
holds for any $x\in\Co$, $t\geq 0$, $q\geq 1$. (The proof is based on
the classical Bernstein inequality for polynomials and 
the polynomial representation of the ${\cal L}$-extremal function
of the disk (see, e.g. [K]).)  Then the estimate of the BMO-norm in
$f|_{\Di(x,t)\cap S}$ follows from inequality (\ref{estbmo}) with $r=1/2$
and $M_{1}-M_{2}=\log 2$.\ \ \ \ \ $\Box$

{\bf Proof of Theorem \ref{reve}.}
As in the proof of Theorem \ref{tbmo} we, first, consider a
local version of the theorem. 
Assume that $K\subset\Di_{r}$ is a compact
from ${\cal A}(d,a,b)$ and $f$, $\Di(x,t)$ satisfy conditions of Lemma \ref{dis}.
Denote $g_{t}=e^{-f'}=e^{f}/e^{f_{x,t}}$. Consider the distribution function
$d_{g}(\lambda):={\cal H}^{d}\{y\in K_{x,t}\ : g_{t}(y)\leq\lambda\}$.
Then from the inequality of Lemma \ref{dis} for $D_{f'}$ we deduce
$$
d_{g}(\lambda)\leq\frac{(4et)^{d}a}{r^{d}d}(\lambda)^{d/(c(M_{1}-M_{2}))}.
$$
Let $g_{*}(s)=\inf\{\lambda\ : d_{g}(\lambda)\geq s\}$. 
From the previous inequality we obtain
$$
g_{*}(s)\geq\left(\frac{sr^{d}d}{(4et)^{d}a}\right)^{c(M_{1}-M_{2})/d}\ .
$$
In particular,
$$
\begin{array}{c}
\displaystyle
\frac{1}{{\cal H}^{d}(K_{x,t})}\int_{K_{x,t}}g_{t}d{\cal H}^{d}=
\frac{1}{{\cal H}^{d}(K_{x,t})}
\int_{0}^{{\cal H}^{d}(K_{x,t})}g_{*}(s)ds\geq\\
\\
\displaystyle
\frac{1}{{\cal H}^{d}(K_{x,t})}\int_{0}^{{\cal H}^{d}(K_{x,t})}\!
\left(\frac{sr^{d}d}{(4et)^{d}a}\right)^{c(M_{1}-M_{2})/d}ds\geq\! 
\frac{1}{1+c(M_{1}-M_{2})/d}
\left(\frac{r^{d}db}{(4e)^{d}a}\right)^{c(M_{1}-M_{2})/d}\ .
\end{array}
$$
Here we used inequality ${\cal H}^{d}(x,t)\geq bt^{d}$. Thus we obtain
\begin{equation}\label{rev1}
\sup_{K_{x,t}}e^{f}\leq (1+c(M_{1}-M_{2})/d)
\left(\frac{(4e)^{d}a}{r^{d}db}\right)^{c(M_{1}-M_{2})/d}
\frac{1}{{\cal H}^{d}(K_{x,t})}\int_{K_{x,t}}e^{f}d{\cal H}^{d}
\end{equation}
which implies the required local reverse H\"{o}lder inequality. In the
general case, we cover again $K$ by a finite number of open disks
$\Di(x_{i},R)$, $i=1,...,N$ such that $f$ is defined in the union of
these disks, the set $\cup_{i=1}^{N}\Di(x_{i},R/2)$ also covers $K$  and
any disk of radius $\leq R/4$ centered at a point of $K$ belongs to one of
$\Di(x_{i},R/2)$.
Then the reverse H\"{o}lder inequality of the form (\ref{rev1}) holds
for any $K_{x,t}=\Di(x,t)\cap K$, $x\in K$, $t\leq R/4$. Assume now that
$t> R/4$ and set
$$
m:=\inf_{x\in K,t>R/4}\left\{
\frac{1}{{\cal H}^{d}(K_{x,t})}\int_{K_{x,t}}e^{f}d{\cal H}^{d}\right\}.
$$
Then $m>0$. Indeed, let
$x_{i},t_{i}>R/4$, be a sequence for which the expression on the right above
converges to $m$. Without loss of generality we may assume also that
$x_{i}$ tends  to $x\in K$ and $t_{i}$ tends to $t\geq R/4$. Then there
is $i_{0}$ such that for any $i\geq i_{0}$, the ball $K_{x_{i},t_{i}}$ contains
$K_{x,R/8}$. Note that $\sup_{K_{x,R/8}}e^{f}>0$ because 
$K_{x,R/8}$ is not a polar set. Then inequality (\ref{rev1}) applied to
$K_{x,R/8}$ and the $d$-regularity of $K$ show that
$$
m\geq
\frac{C}{{\cal H}^{d}(K_{x,R/8})}\int_{K_{x,R/8}}e^{f}d{\cal H}^{d}>0
$$
for a constant $C:=C(K)$. Finally, since $\sup_{K_{x,t}}e^{f}\leq
M:=\sup_{K}e^{f}<\infty$, inequality (\ref{rev1}) for $t>R/4$ is
valid with the constant $M/m$.

The proof of the theorem is complete.\ \ \ \ \ $\Box$
\begin{C}\label{col2}
Assume that a subharmonic function $f$ defined on $\Co$ satisfies
$$
f(z)\leq c'+\log(1+|z|)\ \ \ \ \ (z\in\Co)\ 
$$
for some $c'\in\Re$.
Assume also that $S\in {\cal A}(d,a,b)$. Then for $e^{f}|_{S}$
the reverse H\"{o}lder inequality (\ref{rev1}) holds with the constant
$\frac{c_{1}}{d}\left(\frac{a}{db}\right)^{c_{2}/d}$, where $c_{1}, c_{2}$
are absolute positive constants.
\end{C}
{\bf Proof.}
The proof follows directly from the Bernstein-Walsh inequality (\ref{bernwal})
and Theorem \ref{reve}.\ \ \ \ \ \ $\Box$
\sect{\hspace*{-1em}. Multidimensional Case.}               
In this part we prove the results of section 1.2.
Let $k$ be the number of zeros of $e^{2u}$ in $B_{1/2}$ counting with their
multiplicities (see definition in section 1.2). Here
$u=\frac{1}{2}\log(|f_{1}|^{2}+...+|f_{n}|^{2})$ satisfies inequalities
(\ref{plur}). Below we estimate $k$ by $M,n,r$ only.
\begin{Th}\label{multi}
Given $H>0,d>0$ there exists a system of Euclidean balls
such that
$$
\sum r_{j}^{d}\leq\frac{(2H)^{d}}{d}
$$
where $r_{j}$ are radii of these balls, and
$$
u(z)\geq -M+k\log\frac{H}{e}
$$
everywhere outside these balls in $B_{1/2}$.
\end{Th}
{\bf Proof.}
Let $\xi_{1},...,\xi_{k}$ be zeros of $e^{u}$ in $B_{1/2}$. 
We begin with the following
\begin{Lm}\label{es}
$$
-M+\sum_{i=1}^{k}\log|z-\xi_{i}|\leq u(z)\ \ \ \ \ \ (z\in B_{1/2}).
$$
\end{Lm}
{\bf Proof.}
Without loss of generality we may assume that each of zeros of the system
$F=0$ is of multiplicity 1.
In fact, according to our assumptions image  
$F(B_{1/2})\subset\Co^{n}$ is of complex dimension $n$. In particular,
by Sard's theorem we can 
approximate $F$ by maps $F_{c}=F-c$ where $c$ is a regular value of $F$
close to $0\in\Co^{n}$ and $F_{c}^{-1}(0)$ is a family of zeros of
multiplicity 1. Then we prove
the lemma for $\log|F_{c}|$ and going to the limit as $c\to 0$
obtain the required statement.
Further, observe that $u$ satisfies the complex Monge-Ampere equation
everywhere in $B_{1/2}\setminus\{\xi_{1},...,\xi_{k}\}$. In fact,
$u=\frac{1}{2}F^{*}U$, where $U=\log(\sum_{i=1}^{n}|z_{i}|^{2})$ satisfies 
the Monge-Ampere equation in $\Co^{n}\setminus\{0\}$.
Since $F$ is holomorphic, $u$ satisfies the required equation on 
$B_{1/2}\setminus F^{-1}(0)$. We recall the following result from [BT]:

Assume that $u_{1},u_{2}$ are continuous plurisubharmonic functions in a 
bounded domain $D$ with a compact boundary $K$. Assume also that 
$u_{1}\geq u_{2}$ on $K$ and $u_{1}$ satisfies the complex Monge-Ampere 
equation in  an open neighbourhood of $D$. Then $u_{1}\geq u_{2}$ 
everywhere on $D$.\\
Let $g_{n}=-M+(1+1/n)\sum_{i=1}^{k}\log|z-\xi_{i}|$.
Since by the assumption $\xi_{i}$ is a simple zero of $F$, for any $i$ there
is a ball $B_{r_{n,i}}$ of small radius $r_{n,i}$
centered at $\xi_{i}$ such that 
$g_{n}\leq u$ on its boundary.
Without loss of generality we may assume that
these balls are pairwise disjoint and $r_{n,i}\to 0$ as $n\to\infty$.
Moreover, by definition 
$g_{n}\leq u$ on $S_{1/2}$.
Then according to the above maximal principle,
$g_{n}\leq u$ in $B_{1/2}\setminus (\cup_{i}B_{r_{n,i}})$. It remains to take
the limit as $n\to\infty$ to obtain by continuity
$g\leq u$ in $B_{1/2}$ where $g=-M+\sum_{i=1}^{k}\log|z-\xi_{i}|$.

The lemma is proved.\ \ \ \ \ $\Box$

We now apply Theorem \ref{Cart} to the function $g$ with $X=\Co^{n}$,
$d(x,y)=|x-y|$ and $\mu=\sum_{i=1}^{k}\delta_{\xi_{i}}$. Then 
we obtain 

Given $H>0$, $d>0$, there exists a system of Euclidean balls
such that
$$
\sum r_{j}^{d}\leq\frac{(2H)^{d}}{d}
$$
where $r_{j}$ are radii of these balls, and
$$
g(z)\geq -M+k\log\frac{H}{e}
$$
everywhere outside these balls. Taking into account that $u\geq g$ in
$B_{1/2}$ we obtain the required statement.

The proof of Theorem \ref{multi} is complete.\ \ \ \ \ $\Box$
\begin{R}\label{anyest}
{\rm In the inequality of Theorem \ref{multi} we can take any 
$p\geq k$ instead of $k$. We obtain this replacing inequality of
Lemma \ref{es} by} 
$$
-M+\frac{p}{k}\sum_{i=1}^{k}\log |z-\xi_{i}|\leq u(z)\ \ \ \ \ \
(z\in B_{1/2})
$$
{\rm and then repeating the arguments of the proof of Theorem \ref{multi}
applied to $\frac{k}{p}u$.}
\end{R}
We now estimate the number of zeros $k$.
\begin{Lm}\label{nul}
Under assumptions of Theorem \ref{multi} 
\begin{equation}\label{number}
k\leq c(r,n)e^{(2n-1)M}.
\end{equation}
\end{Lm}
{\bf Proof.}
Let $h=\log(|z_{1}|^{2}+...+|z_{n}|^{2})$. Consider the differential form
$\omega=C(n)(\overline{\partial} h)\wedge
(\overline{\partial}\partial h)^{n-1}$. 
Then we have $d\omega=0$ on $\Co^{n}\setminus\{0\}$ and
for some $C(n)\in\Co$ the Bochner-Martinelli formula is valid
$$
\phi(0)=\int_{\partial D}\phi(\xi)\omega\ .
$$
Here $D$ is a domain containing 0 with a smooth boundary $\partial D$  and
$\phi$ is holomorphic in an open neighbourhood of $\overline{D}$.
Consider now the form $F^{*}\omega$ in $B_{r}$. Since $F$ is a 
holomorphic map, $F^{*}\omega=C(n)(\overline{\partial} F^{*}h)\wedge
(\overline{\partial}\partial F^{*}h)^{n-1}$ and $d(F^{*}\omega)=0$ on
$B_{1/2}\setminus F^{-1}(0)$. In particular, by Stocks' theorem 
$\int_{S_{1/2}}F^{*}\omega=\sum_{i=1}^{k}\int_{S_{i}}F^{*}\omega$, where
$S_{i}$ is a  sphere of a small radius centered at $\xi_{i}$ and
$S_{1/2}=\partial B_{1/2}$. Assume without 
loss of generality that $0$ is a regular value of $F$, i.e., there are
small neighbourhoods of $\xi_{1},...,\xi_{k}$ such that
$F$ maps them biholomorphically to a ball centered at 0.
Assume also that these neighbourhoods contain $S_{1},...,S_{k}$.
Doing in each of these neighbourhoods a holomorphic change of variables and
then applying the Bochner-Martinelli formula we obtain
$\int_{S_{1/2}}F^{*}\omega=k$.
Note that 
$$
F^{*}\omega=C'(n)\frac{F^{*}\sigma}{|F|^{2n}},
$$
where $\sigma=\sum_{i=1}^{n}(-1)^{i-1}\overline{z}_{i}d\overline{z}_{1}
\wedge.._{\hat i}..\wedge\overline{z}_{n}\wedge dz$ and
$dz=dz_{1}\wedge...\wedge dz_{n}$ (see [GH]). From here, Cauchy's inequalities
for the first derivatives of a holomorphic function and the estimate
$|F|\geq e^{-M}$ on $S_{1/2}$ we finally get
$$
k\leq c(r,n)e^{(2n-1)M}.\ \ \ \ \ \Box
$$
\begin{R}\label{unlike}
{\rm Unlike the one-dimensional case, the global Cartan's estimate of Theorem
\ref{multi} do not imply similar local estimates in each ball inside of
$B_{1/2}$ (consider, e.g., function $\log(|z_{1}|^{2}+|z_{2}|^{4})$).
However, under assumptions of Theorem \ref{mcol1} the
multidimensional case is similar to one-dimensional.}
\end{R}

{\bf Proof of Theorem \ref{mcol1}.}
Assume that all zeros of $|F|^{2}$ are elliptic.
\begin{Proposition}\label{m1}
There is a constant $c=c(r,F)>0$ such that for any ball $B(x,t)$ satisfying
$B(x,t)\subset B(x,4r^{2}t)\subset B_{1/2}$ and any $H>0,d>0$ there is a system
of balls such that 
$$
\sum r_{j}^{d}\leq\frac{(8rtH)^{d}}{d},
$$
where $r_{j}$ are radii of these balls, and
$$
\log |F|\geq\sup_{B(x,t)}\log |F|+ k\log\frac{cH}{e}
$$
outside these balls in $B(x,t)$.
\end{Proposition}
{\bf Proof.}

Set
$$
f_{x,t}(z)=\log |F(z)|-\sup_{B(x,4r^{2}t)}\log |F|
$$
and
$$
H_{x,t}=\sup_{t\leq p\leq 2rt}\inf_{z\in S(x,p)}f_{x,t}(z),
$$
where
$S(x,p):=\{z\in\Co^{n}\ :\ |z-x|=p\}$.
Let $K:=\{(x,t)\}$ be the set of centers and radii of balls satisfying
conditions of Proposition \ref{m1}.
\begin{Lm}\label{boun}
$C:=\sup_{(x,t)\in K}H_{x,t}>-\infty $.
\end{Lm}
{\bf Proof.}
Assume that $\{(x_{n},t_{n})\}_{n\geq 1}\subset K$ is a sequence
satisfying
$$
\lim_{n\to\infty}H_{x_{n},t_{n}}=C\ .
$$
Without loss of generality we may assume that $B(x_{n},t_{n})\to
B(x^{*},t^{*})$ in the Hausdorff metric. Further, consider the following
cases.\\
(1) $t^{*}>0$. Then $C=H_{x^{*},t^{*}}$ by continuity. Clearly,
$C>-\infty$ because $F$ has only finite number of zeros in
$B(x^{*},2rt^{*})$.\\
(2) $t^{*}=0$. If $x^{*}$ is not a zero of $F$ then $H_{x^{*},t^{*}}=0$ by
continuity. Assume now that $x^{*}$ is a zero of $F$. Then by ellipticity
of $x^{*}$ the equality
$$
\log |F(x^{*}+t\omega)|=s\log\ t +\log f(\omega)+o(t)\ \ \ \
(0\leq t\leq t_{0},\ \omega\in S^{2N-1}, s\leq k)
$$
holds for a sufficiently small $t_{0}$. Here $0<\inf_{S^{2N-1}}f\leq
\sup_{S^{2N-1}}f<\infty $. From this representation it follows that it suffices
to check the lemma in this case for $|F|=t$, $t\leq t_{0}$, and $x^{*}=0$.
Then a straightforward computation gets
$$
\sup_{B(x_{n},2rt_{n})}\log\ t =\log (|x_{n}|+4r^{2}t_{n})
$$
and so 
\[
H_{x_{n},t_{n}}=
\left\{
\begin{array}{cc}
\displaystyle
\log\frac{|x_{n}|-t_{n}}{|x_{n}|+4r^{2}t_{n}}&
(0\not\in B(x_{n},\frac{(2r+1)t_{n}}{2}));\\
\\
\displaystyle
\log\frac{2rt_{n}-|x_{n}|}{|x_{n}|+4r^{2}t_{n}}&
(0\in B(x_{n},\frac{(2r+1)t_{n}}{2})).
\end{array}
\right.
\]
In the first case $H_{x_{n},t_{n}}$ is a monotonically increasing in $x_{n}$
function. Thus
$$
H_{x_{n},t_{n}}\geq H_{(2r+1)t_{n}/2,t_{n}}=\log\frac{2r-1}{8r^{2}+2r+1}
>-\infty
$$
because $r>1/2$. In the second case
$H_{x_{n},t_{n}}$ is a monotonically decreasing in $x_{n}$ function. This
implies
$$
H_{x_{n},t_{n}}\geq H_{(2r+1)t_{n}/2,t_{n}}=\log\frac{2r-1}{8r^{2}+2r+1}
>-\infty \ .
$$
The proof of the lemma is complete.\ \ \ \ \ $\Box$

We proceed with the proof of Proposition \ref{m1}. According to
Lemma \ref{boun}
there is a sphere $S(x,p)$, $t\leq p\leq 2rt$, such that
$\inf_{S(x,p)}f_{x,t}\geq C>-\infty$. In addition, by conditions of the
theorem $\sup_{B(x,2rp)}f_{x,t}\leq 0$. We set $F'(z)=f_{x,t}(2zp)$,
$|z|\leq r$. Then $F'$ satisfies inequalities (\ref{plur}). Applying
Theorem \ref{multi} to $F'$ and going back to the ball
$B(x,t)\subset B(x,p)$ we obtain

Given $H>0,d>0$, there exists a system of Euclidean balls such that
$$
\sum r_{j}^{d}\leq\frac{(8rtH)^{d}}{d}
$$
where $r_{j}$ are radii of these balls, and
$$
\log |F|\geq \sup_{B(x,4r^{2}t)}\log |F|+C+k\log\frac{H}{e}\geq
\sup_{B(x,t)}\log |F|+k\log\frac{e^{C}H}{e}
$$
outside these balls in $B(x,t)$. We used here that $C\leq 0$.

The proposition is proved.\ \ \ \ \ $\Box$

{\bf Proofs of Theorems \ref{mcol1} and \ref{mcol2}.}
Proofs of these results repeat word-for-word proofs of Theorems
\ref{rem}, \ref{tbmo} and \ref{reve} and might be left to the reader.
\ \ \ \ \ $\Box$

\end{document}